\documentclass[11pt,reqno]{amsart}

\usepackage{amssymb,latexsym,verbatim,pdfsync}

\usepackage{color}

\newcommand{\bbC}{{\mathbb C}}

\newcommand{\bbZ}{{\mathbb Z}}

\def\bR{{\mathbf R}}

\def\cD{{\mathcal D}}
\def\cE{{\mathcal E}}
\def\cH{{\mathcal H}}
\def\cF{{\mathcal F}}

\def\cS{{\mathcal S}}
\def\cU{{\mathcal U}}
\def\cW{{\mathcal W}}

\def\Re{\operatorname{Re}}
\def\Im{\operatorname{Im}}

\def\la{\langle}
\def\ra{\rangle}
\def\raw{\rightarrow}

\def\eps{\varepsilon}
\def\z{\zeta}

\def\ms{\medskip}

\def\endpf{\medskip\hfill $\Box$

\ms

}  
\def\epf{\endpf}

\def\MS{M\"untz--Sz\'asz}

\def\CB{\color{black} }

\newtheorem{thm}{Theorem}[section]
\newtheorem{prop}[thm]{Proposition}
\newtheorem{cor}[thm]{Corollary}
\newtheorem{lem}[thm]{Lemma}

\newtheorem{remark}[thm]{Remark}

\begin{document}

\title[Completeness in Bergman spaces]{Completeness 
on the worm domain and the M\"untz--Sz\'asz
 problem for the Bergman
  space}

\author[S. G. Krantz]{Steven G. Krantz}
\address{Campus Box 1146\\
Washington University in St. Louis\\
St. Louis, Missouri 63130}
\author[M. M. Peloso]{Marco M. Peloso}
\address{Dipartimento di Matematica ``F. Enriques''\\
Universit\`a degli Studi di Milano\\
Via C. Saldini 50\\
I-20133 Milano}
\author[C. Stoppato]{Caterina Stoppato}
\address{Istituto Nazionale di Alta Matematica\\
Unit\`a di Ricerca di Firenze c/o
DiMaI ``U. Dini'' Universit\`a di Firenze\\
Viale Morgagni 67/A\\
I-50134 Firenze\bigskip}

\email{sk@math.wustl.edu\\
marco.peloso@unimi.it\\
stoppato@math.unifi.it}
\thanks{Second author supported in part by the 2010-11 PRIN grant
  \emph{Real and Complex Manifolds: Geometry, Topology and Harmonic Analysis}  
  of the Italian Ministry of Education (MIUR)}
\thanks{Third author supported by the FIRB grant \emph{Differential
    Geometry and Geometric Function Theory} of the MIUR} 
\keywords{Bergman kernel, Bergman projection,
worm domain.}
\subjclass[2010]{32A25, 32A36, 30H20}


\begin{abstract}
In this paper we are concerned with the 
problem of completeness in the Bergman space of the worm domain
$\cW_\mu$ and its truncated version $\cW'_\mu$.  We determine some
orthogonal systems and show that they are not complete, while showing
that the union of two particular of such systems is complete.

In order to prove our completeness result we introduce the {\em
  M\"untz--Sz\'asz} 
for the 1-dimensional Bergman space of the disk $\{ \z:
|\z-1|<1\}$ and find a sufficient condition for its solution.
\end{abstract} 

\maketitle

\section*{Introduction}
\ms

The Diederich--Fornaess worm domain was introduced in
\cite{DFo1} and  is defined for   a given
$\mu>0$ \CB as 
\begin{equation}\label{bdd-worm}
\cW_\mu = \big\{ (z_1,z_2)\in\bbC^2:\,
\big|z_1-e^{i\log|z_2|^2}\big|^2<1-\varphi\big(\log|z_2|^2\big)\big\} \, ,
\end{equation}
where
$\varphi : [-A,A] \to [0,1]$ is a smooth,  convex, even function that
vanishes identically in $[-\mu, \mu]$, with $A>\mu$, $\varphi(A)=1$,
 and increasing on $[\mu, A]$.
As a result, 
$\cW_\mu$ is smooth, pseudoconvex and strictly pseudoconvex
at all points $(z_1,z_2) \in \partial\cW_\mu$ with $z_1\neq
0$.   See \cite{ChSh} for a thorough discussion of 
basic properties of the worm.  The worm turned out to be of fundamental importance in
the theory of geometric analysis in several complex variables, see 
\cite{Ki}, \cite{Barrett-Acta}, \cite{Barrett-Sahutoglu}, \cite{Chr},
\cite{CS}, \cite{KrPe3,KrPe}, 
\cite{KPS1}, \cite{MP1,MP2} and references therein.

For computational purposes, $\cW_\mu$ is often truncated
to the non-smooth bounded domain
\begin{equation}\label{non-smooth-bdd-worm}
\cW'_\mu = \big\{ (z_1,z_2)\in\bbC^2:\,
\big|z_1-e^{i\log|z_2|^2}\big|^2<1,\ \big|\log|z_2|^2\big|<\mu \big\} \,
,
\end{equation}
that is, one replaces the function $\varphi$ with the characteristic
function of the complement of the interval $[-\mu,\mu]$. \ms

In the discussion that follows, we let  $\Omega$ denote  either domain 
$\cW_\mu$ or
$\cW'_\mu$.

In this study we are concerned with the question of finding orthogonal
and complete sets in $A^2(\Omega)$. For details about the notions of 
orthogonality and completeness in a Hilbert space, see Section~\ref{sec:orthogonal}. 
As is well known,  when $\mu\ge\pi$, $\Omega$ has non-trivial {\it Nebenh\"ulle}.
Here the Nebenh\"{u}lle is understood to be the interior of the connected component that contains
$\overline\Omega$ of the intersection of all domains of holomorphy
 containing $\overline\Omega$.
We first show that this easily implies that the closure in
$A^2(\Omega)$  of the holomorphic 
polynomials
is a proper subset of $A^2(\Omega)$.    

Thus we are led to
consider sets of suitable ``monomials'' in $\Omega$ that can be
defined as the holomorphic continuation of non-integral powers
$z_1^\eta$, when $z=(z_1,z_2)$ is initially restricted to
$\Delta\times\{z_2: |z_2|=1\}$, where we denote by
$\Delta$ the disk $\{\z\in\bbC:\, |\z-1|<1\}$.   \ms

We determine some orthogonal sets 
$\{ H_{2k,j}\, , k, j\in\bbZ, k\ge0\}$,
and
$\{ H_{2k+1,j}\, , j\in\bbZ, k\ge0\}$ 
(see Corollary \ref{ort-sys}), and show that 
their union determines a complete set in $A^2(\cW_\mu')$ 
when $\mu>0$   (Theorem
\ref{complete-set}).  We also show that each of the two systems, however, is not
complete (Proposition \ref{incompleteness}).  

In order to prove our completeness result, Theorem \ref{complete-set}, we
prove a result of independent interest, Theorem \ref{M-S}.  We
naturally use the name {\it the \MS\ problem for the Bergman space} for the
question of characterizing the sequences $\{ \lambda_j\}$ in the right
half-plane for which the
sets of powers 
$\{ \z^{\lambda_j-1}\}$ 
form a complete set in
$A^2(\Delta)$.

The classical \MS\ theorem deals with the completeness of
sets of powers $\{t^{\lambda_j-\frac12}\}$ in $L^2\big([0,1]\big)$, 
where again $\lambda_j$ is in the right half plane.  
The solution was provided  by C. M\"untz \cite{Muntz} and by
O. Sz\'asz \cite{Szasz} in two separate papers, where they showed that 
the set $\{t^{\lambda_j-\frac12}\}$  is complete $L^2\big([0,1]\big)$
if and only if $\sum_{j=1}^{+\infty} (1+|\lambda_j|^2)^{-1}
\Re\lambda_j = \infty$ (see also \cite{PW} or \cite{RUD} for a more accessible
reference). 

We find a sufficient condition for the solution of the \MS\ problem
for the Bergman space $A^2(\Delta)$ and use it to prove our completeness
result for $A^2(\cW'_\mu)$.   The \MS\ problem for the Bergman space
has been further studied in \cite{PS}.

Finally we show that the complete set $\{ H_{k,j} , \,
 k, j\in\bbZ, k\ge0\}$,
is not a Schauder basis (Theorem \ref{red-thm}). The definition of Schauder 
basis is recalled in Section~\ref{sec:orthogonal}. 
\ms

\section{Orthogonal sets in $A^2(\cW'_\mu)$}
\label{sec:orthogonal}
\ms

Let  $\mu>0$  and consider the domain $\cW'_\mu$.
The problem we address here is to find a, possibly complete, orthonormal system
for $A^2(\cW'_\mu)$ and consequentely have a way to obtain
an expression for the Bergman kernel.

For the reader's convenience, let us recall a few definitions, for
which we refer to \cite{Young}. A sequence of vectors $\{v_n\}_{n}$ 
in a Banach space $V$ is
\begin{itemize}
\item[{\tiny $\bullet$}] a \emph{Schauder basis} if for each $w\in V$ 
there exists a unique scalar sequence $\{c_n\}_{n}$ such that 
$\sum_{n}c_nv_n$ converges to $w$ with respect to the norm topology. 
\end{itemize}
Now suppose $V$ to be a separable Hilbert space. A sequence 
$\{v_n\}_{n}$ in $V$ is
\begin{itemize}
\item[{\tiny $\bullet$}] an \emph{orthogonal system} if $\la v_m,v_n \ra=0$ 
whenever $m\neq n$;
\item[{\tiny $\bullet$}] an \emph{orthonormal system} if it is orthogonal and 
$\|v_n\|=1$ for all $n$; 
\item[{\tiny $\bullet$}] a \emph{complete system} if $0$ is the only vector in 
$V$ that is orthogonal to $v_n$ for all $n$; or, equivalently, if the linear span 
of $\{v_n\}_{n}$ is dense in $V$.
\end{itemize}
A complete orthonormal system is automatically a Schauder basis and
it is called an \emph{orthornormal basis}. On the other hand, a sequence in 
$V$ that is not orthogonal may be complete without being a Schauder basis; 
in other words, the aforementioned sequences $\{c_n\}_{n}$ may exist for all 
$w\in V$ without being unique.   \medskip

In the analysis on the worm domains $\cW_\mu$ and $\cW_\mu'$ a special
role is played by the functions
\begin{equation}\label{def-eta}
E_\eta(z) = e^{\eta L(z)}   \, ,
\end{equation}
where
\begin{equation}\label{def-L}
L(z) = \log\big(z_1e^{-i\log|z_2|^2}\big)+i\log|z_2|^2\, ,
\end{equation}
and $\log$ denotes the principal branch of the logarithm, so that
$$
E_\eta(z_1,z_2) = \big(z_1e^{-i\log|z_2|^2}\big)^\eta
e^{i\eta\log|z_2|^2}\, .
$$

The function $L$ is well defined and holomorphic in a domain
containing $\cup_\mu \cW'_\mu$ (see \cite{KPS1}, Lemma 1.2
and Proposition 1.3).  Moreover, we point out that the fiber of $\cW'_\mu$ over each $z_1 \in
D(0,2)\setminus\{0\}$ is not connected and that $L(z)$ is locally constant
in $z_2$, but not constant. The same happens with $E_\eta(z) $ for $\eta
\in \bbC \setminus \bbZ$, while $E_k(z) = z_1^k$ for all $k \in \bbZ$, $z \in
\cW'_\mu$.   Hence the functions $E_\eta$ are the analytic
continuation to $\cW'_\mu$ of the monomial $z_1^\eta$ defined in
$\cW'_{\pi/2}$ using the principal branch of the logarithm.\ms

It is well known that the functions that are holomorphic in a
neighborhood of the closure $\overline{\cW_\mu }$ are not dense in
$A^2(\cW_\mu)$.  Since a proof of this fact does not explicitly appear in
the literature, we prove the following result that applies to both
domains $\cW_\mu$ and $\cW_\mu'$.

\begin{prop}\label{non-density}	  \sl
 Let  $\mu\ge2\pi$ and let  
$A^2(\overline{\cW_\mu})$ denote the closure in
  $A^2(\cW_\mu)$ of the functions that are holomorphic in a neighborhood of
  $\cW_\mu$.  Then, if $f\in A^2(\overline{\cW_\mu})$, then $f$ is
  holomorphic on 
$\widehat{\cW_\mu}$, where
\begin{equation}\label{w-hat}
\widehat{\cW_\mu} 
\supseteq \bigcup_{-\mu\le a\le \mu-2\pi}  
\big\{ (z_1,z_2): a<\log |z_2|^2 <a+2\pi,\ |z_1-e^{ia} |<1
\big\}. 
\end{equation}
Therefore $A^2(\overline{\cW_\mu})\subsetneqq A^2(\cW_\mu)$. 

The same conclusions hold true with $\cW_\mu'$ in place of $\cW_\mu$.
\end{prop}

In particular, the polynomials are not dense in either $A^2(\cW_\mu)$
or  $A^2(\cW_\mu')$.    By contrast, D. Catlin \cite{Cat3} showed that
for every smoothly bounded pseudoconvex domain $\Omega$,
the holomorphic functions in $C^\infty(\overline{\Omega})$ are dense
in $A^2(\Omega)$.  
 
\proof
It suffices to prove the result in the case of $\cW_\mu'$, since the
same argument can be repeated {\it verbatim} for $\cW_\mu$.

Suppose $f$ is holomorphic in a neighborhood of $\overline{\cW_\mu'}$ 
and let 
$$
\cE_a =\big\{ (0,z_2):\, a \leq \log|z_2|^2 \leq a+2\pi \big\} \bigcup
 \big\{ (z_1,z_2):\, \log|z_2|^2 =a \ \hbox{or} \ a+2\pi,\ |z_1-e^{ia}| \le 1\big\}
\,.
$$
Set
$$
F_a(z_1,z_2) = 
\frac{1}{2\pi i} \int_\gamma \frac{f(z_1,\z)}{\z-z_2}\, d\z \,,
$$
where $\gamma$
 is the oriented boundary of the annulus $\{z_2 \in \bbC :
 a<\log|z_2|^2<a+2\pi\}$. 
Then $F_a$ is holomorphic on the set
$$
\widetilde \cE_a=
\big\{ (z_1,z_2):\, a< \log|z_2|^2<a+2\pi,\ |z_1-e^{ia}| <1
\big\}\, .
$$
However, $F_a(0,z_2)= f(0,z_2)$, since $f$ is holomorphic in a
neighborhood of $\{(0,z_2): \big| \log|z_2|^2\big| \le \mu\}$ and
$-\mu\le a\le \mu-2\pi$ implies $a+2\pi \leq \mu$.  It follows that $F_a$ is
a holomorphic extension of $f$ to the set $\widetilde \cE_a$ and thus
$f$ extends holomorphically to an open set $\widehat{\cW_\mu'}$
containing the right-hand side of \eqref{w-hat}.

Now suppose that $\{f_n\}$ are holomorphic in a a neighborhood of 
$\overline{\cW_\mu'}$ and that $f_n\to f$ in $A^2(\cW_\mu')$.  Then
$f_n\to f$ uniformly on the compact sets $\{ \log |z_2|^2=a,\,
a+2\pi;\ |z_1-e^{ia}|\le 1-\delta\}$, for $-\mu< a<  \mu-2\pi$.  By
Cauchy's formula, $\{f_n\}$ is Cauchy in uniform norm also on the sets
$$
\big\{ a\le \log |z_2|^2\le a+2\pi,\ |z_1-e^{ia}|\le 1-\delta
\big\} \, .
 $$
Therefore $f$ extends holomorphically to the set
on the right-hand side of \eqref{w-hat}.

Finally, the functions $E_\eta$ with $\eta$ not an integer cannot be
extended  holomorphically to any of the sets $\widetilde \cE_a$,
  so that
$A^2(\overline{\cW_\mu'})\subsetneqq A^2(\cW_\mu')$.  

It is immediate to check that the arguments above apply to the case of
$\cW_\mu$ as well.
\qed
\ms

Thus, we are led to consider the set of ``monomials'' of the form
$\{E_{\eta_j}(z)z_2^j\}_{j\in\bbZ}$ and ask whether these are
orthogonal, and/or complete, for some choice of values
$\eta_j\in\bbC$.  
\ms

We denote by  $dA$ the Lebesgue
measure in the complex plane. 

\begin{lem}\label{Beta-fnc}  \sl
Let $\Re\alpha,\Re\beta>-1$.  Then
$$
\int_\Delta \z^\alpha\overline{\z^\beta}\, dA(\z)
= \pi \frac{\Gamma(\alpha+\overline\beta+2)}{
\Gamma(\alpha+2)\Gamma(\overline\beta+2)}\, .
$$
In particular, $\z^\alpha$ and $\z^\beta$ are never orthogonal to each
other in $A^2(\Delta)$.
\end{lem}
\proof
We have
\begin{align*}
\int_\Delta \z^\alpha\overline{\z^\beta}\, dA(\z)
& = \int_0^2 \int_{-\cos^{-1}(r/2)}^{\cos^{-1}(r/2)}
\big(re^{i\theta}\big)^\alpha \big(re^{-i\theta}\big)^{\overline\beta}
\, d\theta\, r\, dr\\
& = \int_0^2 r^{\alpha+\overline\beta+1} \int_{-\cos^{-1}(r/2)}^{\cos^{-1}(r/2)}
e^{i\theta(\alpha-\overline\beta)} \, d\theta dr\\
& = \frac{2}{i(\alpha-\overline\beta)}
\int_0^2 r^{\alpha+\overline\beta+1} \sinh \big(
\cos^{-1}(r/2)i(\alpha-\overline\beta)\big) \, dr\\
& = \frac{4}{i(\alpha-\overline\beta)}
\int_0^{\pi/2} \big(2\cos s\big)^{\alpha+\overline\beta+1} 
\sin s \sinh \big(is(\alpha-\overline\beta)\big)\, ds \\
& = \frac{2^{\alpha+\overline\beta+3}}{\alpha+\overline\beta+2}
\int_0^{\pi/2} \big(\cos s\big)^{\alpha+\overline\beta+2} 
\cosh \big(is(\alpha-\overline\beta)\big)\, ds \\
& = \frac{2^{\alpha+\overline\beta+3}}{\alpha+\overline\beta+2}
\int_0^{\pi/2} \big(\cos s\big)^{\alpha+\overline\beta+2} 
\cos \big(s(\alpha-\overline\beta)\big)\, ds \, .
\end{align*}

Now we use
\cite[(9) p. 391]{tables} and, denoting by $B$ the beta function, we obtain that
\begin{align*}
\int_\Delta \z^\alpha\overline{\z^\beta}\, dA(\z)
& = \frac{\pi}{(\alpha+\overline\beta+2)(\alpha+\overline\beta+3)}
\cdot\frac{1}{B(\alpha+2,\overline\beta+2)}\\
& = \pi \frac{\Gamma(\alpha+\overline\beta+2)}{
\Gamma(\alpha+2)\Gamma(\overline\beta+2)}\, ,
\end{align*}
as we wished to prove. \ms
\epf

For a given bounded  domain $\Omega$ in $\bbC^2$ that is rotationally
invariant\footnote{These are called {\it Hartogs} domains.} in
the second variable $z_2$, such as $\cW_\mu$ and $\cW_\mu'$,  using the 
Fourier expansion in $z_2$, the Bergman space
$A^2(\Omega)$ decomposes as an orthogonal sum
\begin{equation}\label{oplus-sum}
A^2(\Omega) =\bigoplus_{j\in\bbZ} 
\cH^j\, .
\end{equation}
Here 
$$
\cH^j =\bigl\{ F\in A^2(\Omega):\, 
F(z_1,z_2) =f(z_1,|z_2|)z_2^j\bigr\} \, ,
$$
where $f$ is holomorphic in $z_1$ and locally constant in $|z_2|$. 
The orthogonal
projection of $A^2$ onto $\cH^j$ is given by
$$
Q_jF(z_1,z_2) = \frac{1}{2\pi} \int_{-\pi}^\pi
F(z_1,e^{it}z_2)e^{-ijt}\, dt\, .
$$
Then we set
$$
f(z_1,|z_2|) = \frac{Q_j F(z_1,z_2)}{z_2^j}\, ,
$$
and observe that the right-hand side is 
holomorphic in $\Omega$, but depends only on the modulus of $z_2$.
Hence $f$ is locally constant in $|z_2|$.  In general, it will be
constant only if the fibers over a point $(z_1,z_2)$ in $\Omega$ with
$z_1$ fixed, is a connected set in the $z_2$-plane.
\ms

Let us go back to the case of $\cW'_\mu$ and
let
 $F,G\in \cH^j$, $F(z_1,z_2) =f(z_1,|z_2|)z_2^j$ and $G(z_1,z_2)
=g(z_1,|z_2|)z_2^j$.  We have
\begin{align}
\la F,\, G\ra_{A^2(\cW'_\mu)} 
& =  \int_{|\log|z_2|^2|<\mu} 
\int_{|z_1-e^{i\log|z_2|^2}|<1} f(z_1,|z_2|)
\overline{g(z_1,|z_2|)} |z_2|^{2j}\,  dA(z_1)dA(z_2) \notag\\
& = 2\pi \int_{|\log r^2|<\mu} 
\int_{|z_1-e^{i\log r^2}|<1} f(z_1,r)
\overline{g(z_1,r)} \,  dA(z_1)\, r^{2j+1} dr  \notag\\
& = \pi \int_{|s|<\mu} 
\int_{|z_1-e^{is}|<1} f(z_1,e^{s/2})
\overline{g(z_1,e^{s/2})} \,  dA(z_1)\, e^{s(j+1)} ds \notag\\
& = \pi\int_\Delta \int_{|s|<\mu} 
 f(\z e^{is},e^{s/2})
\overline{g(\z e^{is},e^{s/2})}\, e^{s(j+1)} ds \,  dA(\z)\, .\label{F-G}
\end{align}

Let 
\begin{equation}\label{nu}
\nu=\pi/2\mu\, .
\end{equation}
be the reciprocal of the winding number of $\cW'_\mu$, as
also defined in \cite{Barrett-Acta} and
let $h$ be the entire function
\begin{equation*}
h(z)=\frac{\sinh
  [\mu(j+1+iz)]}{j+1+iz} \, .
\end{equation*}
Define
\begin{equation}\label{gamma-alpha-beta}
\gamma_{\alpha\beta}
= h(\alpha-\overline\beta)
\, .
\end{equation}

\begin{prop}\label{prop1}  \sl
 Let $\mu>0$.  
For $\alpha\in\bbC$ and $j\in\bbZ$ let 
$F_{\alpha,j}(z_1,z_2)=E_\alpha(z) z_2^j$.  Then $F_{\alpha,j}\in A^2(\cW'_\mu)$ if
and only if $\Re\alpha>-1$.  Moreover, if $\Re\alpha,\Re\beta>-1$
then
$$
\la F_{\alpha,j},\, F_{\beta,j}\ra_{A^2(\cW'_\mu)}
= (2\pi)^2 \gamma_{\alpha\beta} \frac{\Gamma(\alpha+\overline\beta+2)}{
\Gamma(\alpha+2)\Gamma(\overline\beta+2)}\, .
$$
In particular, 
$\la F_{\alpha,j},\, F_{\beta,j}\ra_{A^2(\cW'_\mu)} =0$
if and only if
\begin{equation}\label{orth-cond}
\alpha-\overline\beta= 2k\nu +i(j+1)\qquad
\text{with\ } k\in\bbZ\setminus\{0\}\, .
\end{equation}
\end{prop}
\proof
We compute $\la F_{\alpha,j},F_{\beta,j}\ra_{A^2(\cW'_\mu)}$. 
Starting from \eqref{F-G} we obtain
\begin{align}
\la F_{\alpha,j},\, F_{\beta,j}\ra_{A^2(\cW'_\mu)}
& = \pi\int_\Delta \int_{|s|<\mu} 
 E_\alpha(\z e^{is}, e^{s/2})
\overline{E_\beta(\z e^{is}, e^{s/2})}\, e^{s(j+1)} ds \,  dA(\z)\notag\\
& = \pi\int_\Delta \z^\alpha\overline{\zeta^\beta} \int_{|s|<\mu} 
 e^{is(\alpha-\overline\beta)}\, e^{s(j+1)} ds \,  dA(\z)\notag\\
& = 2\pi\gamma_{\alpha\beta} \int_\Delta
\z^\alpha\overline{\z^\beta} \, dA(\z)\, ,
\label{F-G2}
\end{align}
where 
\begin{align*}
\gamma_{\alpha\beta}
& := \textstyle{\frac12} 
\displaystyle{\int_{|s|<\mu}} 
 e^{s(j+1+i(\alpha-\overline\beta))}\,  ds \\
& = \begin{cases}
\frac{\sinh
  (\mu(j+1+i(\alpha-\overline\beta)))}{j+1+i(\alpha-\overline\beta)}
& \text{if}\ j+1+i(\alpha-\overline\beta)\neq 0
\cr 
\mu & \text{if}\ j+1+i(\alpha-\overline\beta)=0\, ,
\end{cases}
\end{align*}
as claimed.

Therefore,
$\gamma_{\alpha\beta}=0$ if and only if
$$
\mu(j+1+i(\alpha-\overline\beta))=k\pi i\qquad\qquad
\text{for\ }k\in\bbZ\setminus\{0\}\, ,
$$
that is, 
\begin{equation}\label{gamma=0}
\alpha-\overline\beta= 2k\nu +i(j+1) 
\qquad\qquad\text{for\ }k\in\bbZ\setminus\{0\}\, .\ms
\end{equation}

Notice that, when $\alpha=\beta$,   the previous computation gives
\begin{align*}
\| E_\alpha(z) z_2^j\|_{A^2(\cW'_\mu)}^2
& = 2\pi \gamma_{\alpha,\alpha} \int_\Delta
|\z^\alpha|^2 \, dA(\z) \\
& =  2\pi \gamma_{\alpha,\alpha} \int_{|\z|<1}
e^{2[\Re\alpha\log|\z-1|-\Im\alpha \arg(\z-1)]}  \, dA(\z) \\
& = 2\pi \gamma_{\alpha,\alpha} \int_{|\z|<1}|\z-1|^{2\Re\alpha}
e^{-2\Im\alpha \arg(\z-1)}  \, dA(\z) \, ,
\end{align*}
which is finite if and only if $\Re\alpha>-1$.  This proves the first part of the
statement. The second part now follows from Lemma \ref{Beta-fnc}.
\epf

The following corollaries now follow at once.

\begin{cor}\label{norm-F-alpha}	  \sl
  Then, for $\Re\alpha>-1$, we have
that  
$F_{\alpha,j}\in A^2(\cW'_\mu)$ and  
$$
\|F_{\alpha,j}\|^2_{A^2 (\cW'_\mu)  }
= (2\pi)^2 \frac{\sinh[\mu(j+1-2\Im\alpha)]}{j+1-2\Im\alpha}\,
\frac{\Gamma(2+2\Re\alpha)}{|\Gamma(2+\alpha)|^2}\, .
$$
\end{cor}

For $c_0>-1$, and $\ell=0,1,2,\dots$  we set 
\begin{equation}\label{H-ell-j-def}
H _{\ell,j}(z_1,z_2) =
E_{c_0+\nu\ell+i(j+1)/2}(z)z_2^j\, .
\end{equation}

\begin{cor}\label{ort-sys}  \sl
For $\mu>0$, each of the two sets
\begin{equation}\label{two-sets}
\big\{  H_{2k,j}\, ,\, j\in\bbZ,\,  k=0,1,2,\dots\big\},
\quad\text{and}\quad
\big\{  H_{2k+1,j}\, ,\, j\in\bbZ,\,  k=0,1,2,\dots\big\},
\end{equation}
is an orthogonal system in $A^2(\cW'_\mu )$.
\ms
\end{cor}

\section{The \MS\ problem for the Bergman
  space}
\ms

In endeavoring to establish whether the system $\{ F_{\alpha,j}\}$ is
complete we are led to consider the 
\MS\  problem for the Bergman space.
\ms

Recall that we set $\Delta=\{\z:\, |\z-1|<1\}$.  We consider a set of
functions $\{\z^{\lambda_k}\}$, $k=1,2,\dots$ and would like to find a
necessary and sufficient condition for this set to
be a {\it complete set} in $A^2(\Delta)$, that is, its linear span to be dense
in $A^2(\Delta)$.

\begin{thm}\label{M-S}	 \sl
Let $\cS$ be the subset of $A^2(\Delta)$ whose elements are the
functions $\zeta^{\lambda_k}$ for $k = 0,1,2,\dots$, 
where $\lambda_k=ak+c_0+ib$, $0<a<1$, $c_0>-1$ and $b\in\bR$.
Then $\cS$ is a  complete set in  $A^2(\Delta)$. 
\end{thm}

\proof
Consider the biholomorphic map $C:\cU\raw \Delta$ given by
$$
C(\omega)= \frac{2i}{i+\omega}
$$
of the upper half plane $\cU$ onto $\Delta$.  Then
$$
T: A^2(\Delta)\ni f \mapsto (f\circ C)C' \in A^2(\cU)
$$
is a surjective isometry.  Next,  by the
Paley--Wiener theorem, the  Fourier transform $\cF$, given by 
$$
(\cF g) (\xi) = \frac{1}{2\pi}\int_{-\infty}^{+\infty} g(t)e^{-i\xi t}\, dt \, ,
$$
provides a surjective isometry of $A^2(\cU)$ onto 
$L^2\big( (0,+\infty),d\xi/\xi\big)$, see e.g. \cite{5p}.
Therefore, $ \{\z^{\lambda_k}\}$ will be complete in
$A^2(\Delta)$ if and only if $\big\{\cF \big(T
\z^{\lambda_k}\big)\big\}$ is complete in $L^2\big(
  (0,+\infty),d\xi/\xi\big)$.

Now, 
$$
T(\z^\lambda)(\omega) =-
\frac{(2i)^{\lambda+1}}{(i+\omega)^{\lambda+2}}
\, , 
$$ 
while
\begin{align*}
\cF \Bigl( (u+i)^{-(\lambda+2)} \Bigr) (\xi)
& = \frac{1}{2\pi} \int_{-\infty}^{+\infty} 
\frac{1}{(u+i)^{\lambda+2}}
e^{-i\xi u} \, du \\
& = \frac{1}{i^{\lambda+2}\Gamma(\lambda+2)}
\xi^{\lambda+1}e^{-\xi}\chi_{(0,+\infty)}(\xi) \, ,
\end{align*}
(see also \cite{DGM}, Lemma 1).

Therefore,
$$
\cF\big(T(\z^\lambda)\big) (\xi) 
= -\frac{(2\xi)^{\lambda+1}e^{-\xi}}{i\Gamma(\lambda+2)}
\chi_{(0,+\infty)}(\xi) \, .
$$
Hence the set $\{\z^{\lambda_k}\}$ is complete in $A^2(\Delta)$ if and
only if the set $\{\xi^{\lambda_k+1}e^{-\xi}\}$ is complete in 
$L^2\big(
  (0,+\infty),d\xi/\xi\big)$, that is, the set 
 $\{\xi^{\lambda_k+\frac12}e^{-\xi/2}\}$ is complete in 
$L^2\big((0,+\infty) \big)$.  

Next, we consider the transformation $\xi\mapsto \xi^a=t$ of
$(0,+\infty)$ onto itself and the induced isometry $\Lambda$ of $L^2\big(
  (0,+\infty)\big)$ onto itself given by
$$
(\Lambda \psi) (t) = \textstyle{\sqrt{\frac1a}} \psi (t^{1/a})
\Bigl(t^{\frac12(\frac1a - 1)}
\Bigr)\, .
$$
Under such a transformation, since $\lambda_k=ak+c_0 +ib$, we see that 
$\{\xi^{\lambda_k+\frac12}e^{-\xi/2}\}$ is complete in 
$L^2\big((0,+\infty) \big)$ if and only if
$\{t^k t^\alpha e^{-\frac12 t^{1/a}}\}$ is complete in 
$L^2\big((0,+\infty) \big)$, where 
$$\alpha=
\frac{c_0+1}{a}-\frac12+i\frac{b}{2a}\, .
$$

We know
from \cite[Thm. 5.7.1]{Szego} that the system $\{t^{n+c}e^{-t/2}:\, n=0,1,2,\dots\}$,
with $c>-1/2$, 
is complete in $L^2\big((0,+\infty) \big)$.  
Thus, if $\psi\in
L^2\big((0,+\infty) \big)$ is orthogonal to $t^k t^\alpha
e^{-\frac12 t^{1/a}}$ for all $k=0,1,2,\dots$ it follows that 
$$
\int_0^{+\infty} t^{k+\Re \alpha}e^{-t/2} 
\overline{t^{-i\Im\alpha}e^{\frac12(t-t^{1/a})}\psi(t)}\, dt =0
$$
for $k=0,1,2,\dots$.  Since $0<a<1$, $e^{\frac12(t-t^{1/a})}$ is bounded and
also $\Re\alpha>-\frac12$, we obtain that 
$t^{-i\Im\alpha}e^{\frac12(t-t^{1/a})}\psi=0$, that is, $\psi=0$.
This concludes the proof.
\ms
\epf

\section{Complete  sets in
  $A^2(\cW'_\mu)$}\label{comp-sets-section} 
\ms

From our \MS\ Theorem \ref{M-S} for the Bergman space
$A^2(\Delta)$, 
we obtain the following density 
result in $A^2(\cW'_\mu)$.
\begin{thm}\label{complete-set}	 \sl 
 Let $\mu>\pi/2$.  
Let $ H _{\ell,j}(z_1,z_2)$ be as in \eqref{H-ell-j-def}. 
Then
$\{ H _{\ell,j}\}_{\ell, j\in\bbZ,\,\ell\ge0}$, 
is a complete set in $A^2(\cW_\mu')$.
\end{thm}
Notice that the set $\{ H _{\ell,j}\}$, $\ell, j\in\bbZ$, $\ell\ge0$, is the
union of the two sets in \eqref{two-sets}.

\proof  We wish to show
that if $F\in A^2(\cW_\mu')$ is orthogonal to $ H _{\ell,j}$, for $\ell,
j\in\bbZ$, $\ell\ge0$, then $F$ is identically zero.  It suffices to
show that, for each $j\in\bbZ$ fixed, any function $F\in \cH^j$
orthogonal to $ H _{\ell,j}$ for all $\ell\ge0$, is identically zero.

Writing $F(z_1,z_2)=f(z_1,|z_2|)z_2^j$, from \eqref{F-G2} we then have
\begin{align}
0 & = \la F,\,  H _{\ell,j}\ra_{A^2(\cW'_\mu)} =  \pi\int_\Delta \int_{|s|<\mu} 
f(\z e^{is},e^{s/2})
\overline{E_{c_0+\nu\ell+i(j+1)/2}(\z e^{is}, e^{s/2})}\, e^{s(j+1)} ds \,
dA(\z)\notag \\
& = \pi\int_\Delta \int_{|s|<\mu} 
f(\z e^{is},e^{s/2}) \,
 e^{s[(j+1)/2 +i(c_0+\nu\ell)]} \,  ds \, \overline{\z^{c_0+\nu\ell+i(j+1)/2}} 
dA(\z)\, ,\label{orthogonality-cond}
\end{align}
for $\ell=0,1,\dots\,$.   
Notice that the function
\begin{align}
Tf(\z,w)
&
=
\int_{|s|<\mu} f(\z e^{is},e^{s/2}) e^{s[(j+1)/2 +ic_0]}\,
 e^{iw} \,  ds \notag\\
&  =\cF\big( f(\z e^{is},e^{s/2}) e^{s[(j+1)/2
   +ic_0]}\chi_{\{|s|<\mu\}} \big)(w) 
\end{align}
is analytic in $\z\in\Delta$, and by the
Paley--Wiener theorem \cite{PW}, is an entire function in $w$ of
exponential type at most $\mu$.  
Moreover, 
the function
$$
w\mapsto \pi\int_\Delta Tf(\z,w) \overline{\z^{c_0+\nu\ell+i(j+1)/2}} 
dA(\z)
$$ 
is again an entire function of
exponential type at most $\mu$ and by \eqref{orthogonality-cond} 
 it vanishes at the points
$w_\ell = \nu\ell$.  
Observe that
\begin{equation}\label{Fuchs-cond}
\limsup_{r\to+\infty} \frac{ \exp\{ 2 \sum_{\nu\ell<r} 1/(\nu\ell) \}
}{r^{2\mu/\pi}} 
=  \Big(  \frac{1}{\nu} \Big)^{\frac2\nu} 
\lim_{r\to+\infty} r^{1/\nu} = +\infty\, .
\end{equation}
By a classical result of Fuchs \cite{Fuchs}, we know that an entire
function of type $\mu$ whose zero set $\{w_\ell=\nu\ell\}$ satisfies
\eqref{Fuchs-cond} must vanish identically, that is,
$$
\int_\Delta Tf(\z,w) \overline{\z^{c_0+\nu\ell+i(j+1)/2}} 
dA(\z) =0\, ,
$$ 
for $\ell=0,1,\dots$.  
Since $\mu>\pi/2$ we have $\nu<1$ and by Theorem
\ref{M-S} it now follows that
$Tf(\cdot,w)=0$, hence, 
\begin{equation}\label{reduction}
\int_{|s|<\mu} 
f(\z e^{is},e^{s/2}) \,
 e^{s[(j+1)/2 +ic_0]} e^{isw} \,  ds =0 \, ,
\end{equation}
for all $\z\in\Delta$ and $w\in\bbC$. 
  This implies that $f$ vanishes
identically and we are done. \ms\qed

Notice that, had we considered either of the orthogonal systems
	mentioned in Corollary \ref{ort-sys}, we would have ended up with
	the points $\{w_{2\ell}\}$ only, or with $\{w_{2\ell+1}\}$. The analog of 
	Condition \eqref{Fuchs-cond} would not have been satisfied and 
	we could not have proved completeness using this approach.
	In fact, we are going to show in the next proposition that each of the 
	two systems is incomplete.\ms

It is also worth mentioning that the worm domains $\cW'_\mu$ are
increasingly badly behaved as $\mu$ becomes large.  On the other hand,
the proof of 
our density result breaks down when $\mu\le\pi/2$.
This is somewhat surprising, 
 since when $\mu\le\pi/2$ the fibers over $z_1$ are connected,
the geometry of the domain is much simpler and in principle it should
be easier to obtain such results on $\cW'_\mu$ when $\mu\le\pi/2$.

\ms 

\begin{prop}\label{incompleteness}   \sl
Let  
$\{ H_{\ell,j}\}_{\ell, j\in\bbZ,\,\ell\ge0}$ be as in Theorem
\ref{complete-set}.  
Then, for each $m$ fixed,
\begin{equation}\label{strict-ineq}
 \|  H_{2m+1,j} \|^2_{A^2(\cW'_\mu)}  > 
\sum_{j'=-\infty}^{+\infty} \sum_{k=0}^{+\infty}  \frac{1}{\|  
  H_{2k,j'} \|^2_{A^2(\cW'_\mu)}}
\big| \big\la  H_{2m+1,j},\, H_{2k,j'}  \big\ra\big|^2\, ,
\end{equation}
and, analogously, for each $k$ fixed,
$$
\| H_{2k,j}\|^2_{A^2(\cW'_\mu)} > 
\sum_{j'=-\infty}^{+\infty} \sum_{m=0}^{+\infty}  \frac{1}{\| H_{2m+1,j'}\|^2_{A^2(\cW'_\mu)}}
\big| \big\la H_{2k,j},\, H_{2m+1,j'}\big\ra\big|^2\, .
$$

Hence, neither system $\{H_{2k,j}\}$ nor $\{H_{2m+1,j}\}$ is complete 
in $A^2(\cW'_\mu)$. 
\end{prop} 

\proof 
By orthogonality, it suffices to consider the case $j'=j$ and,  
 dropping the index $j$,  we  write $F_k=H_{2k,j}$ and
$G_m=H_{2m+1,j}$. 
By Proposition \ref{prop1} we have that
\begin{align*}
\la G_m,\, F_k \ra_{A^2(\cW'_\mu)} 
 & = (2\pi)^2 \frac{\sin \big[ \mu\big(2k-(2m+1)\big)\nu\big]}{
  \big(2k-(2m+1)\big)\nu} \\
& \qquad\qquad \times \frac{ 
\Gamma\big(2c_0+2 +  \big(2(k+m)+1\big)\nu\big)  }
{\Gamma\big(c_0+2+ 2k\nu +i\textstyle{ \frac{j+1}{2}}  \big)
\Gamma\big(c_0+2+ (2m+1)\nu -i\textstyle{ \frac{j+1}{2}}  \big)} \, , \\
\| F_k \|_{A^2(\cW'_\mu)}^2
& = (2\pi)^2 \mu \frac{
\Gamma\big(2c_0+2  + 4k\nu\big)  }
{\big| \Gamma\big(c_0+2+ 2k\nu +i\textstyle{ \frac{j+1}{2}}  \big)
  \big|^2} \, ,\\
\| G_m \|_{A^2(\cW'_\mu)}^2
& = (2\pi)^2 \mu \frac{
\Gamma\big(2c_0+2  +2(2m+1)\nu\big)  }
{\big| \Gamma\big(c_0+2+ (2m+1)\nu -i\textstyle{ \frac{j+1}{2}}  \big)
  \big|^2}\, .
\end{align*}

Therefore, \eqref{strict-ineq} is equivalent to
\begin{align}
& \mu 
\Gamma\big(2c_0+2  +2(2m+1)\nu\big) \notag \\ 
 & >  \sum_{k=0}^{+\infty} \bigg[ \frac{
\sin \big[ \mu\big(2k-(2m+1)\big)\nu\big]}{
  \big(2k-(2m+1)\big)\nu} \bigg]^2 
\frac{  \Gamma\big(2c_0+2 +  \big(2(k+m)+1\big)\nu\big) ^2 
}{ \mu \Gamma\big(2c_0+2   +4k\nu\big) } \, ,   \label{strict-ineq-2} 
\end{align}
 which in turn is implied by, 
\begin{align}
 1
& > \sum_{k=0}^{+\infty} \frac{1}{
  \big[ \mu \big(2k-(2m+1)\big)\nu\big]^2}  
\frac{ \Gamma\big(2c_0+2 +  \big(2(k+m)+1\big)\nu\big) ^2 
}{  \Gamma\big(2c_0+2  +2(2m+1)\nu\big) 
\Gamma\big(2c_0+2   +4k\nu\big) } \notag \\
& = \sum_{k=0}^{+\infty} \frac{1}{ \pi^2 \big(k-\frac{2m+1}{2}\big)^2}  
\frac{ \Gamma\big(2c_0+2  + \big(2(k+m)+1\big)\nu\big) ^2 
}{  \Gamma\big(2c_0+2  +2(2m+1)\nu\big) 
\Gamma\big(2c_0+2   +4k\nu\big) }
\, . \label{strict-ineq-3}
\end{align}

Now, on the one hand the right-hand side in \eqref{strict-ineq-3} 
is less than or equal to
$$
\sum_{k=0}^{+\infty} \frac{1}{ \pi^2 \big(k-\frac{2m+1}{2}\big)^2}  \, ,
$$
since for all $x,y>0$, $c\ge0$, 
\begin{align*}
\Gamma(c+x+y)^2 
& =\left(\int_0^{+\infty} t^{x+y+c-1} e^{-t} dt \right)^2  \\
&  \leq \left(\int_0^{+\infty} t^{2x+c-1} e^{-t} dt\right)
\left(\int_0^{+\infty} t^{2y+c-1} e^{-t} dt \right)
= \Gamma(c+2x) \Gamma(c+2y)\, .
\end{align*}

On the other hand, we claim that
$$
\pi^2  > \sum_{k=0}^{+\infty} \frac{1}{\big(k-\frac{2m+1}{2}\big)^2}.
$$
Indeed, setting $h(w) = \pi \cot(\pi w)$ and $Q(w) =
\big(w-\frac{2m+1}{2} \big)^2$, we have
$$
\operatorname{Res}\left(\frac{h}{Q},k\right) 
= \frac{1}{Q(k)} = \frac{1}{\big(k-\frac{2m+1}{2}\big)^2}
$$
for all $k \in \bbZ$ and
$$
\operatorname{Res}\left(\frac{h}{Q},\frac{2m+1}{2}\right) 
= \lim_{w \to \frac{2m+1}{2} } \frac{\pi}{\sin(\pi w)} 
\frac{\cos(\pi w)}{w-\frac{2m+1}{2}} 
= -\pi^2.
$$
The fact that
$$
0 = \lim_{n \to +\infty} \int_{\partial D(0,n+\frac{1}{2})}
\frac{h(w)}{Q(w)} dw 
= 2\pi i \lim_{n \to +\infty}
\left[\operatorname{Res}\left(\frac{h}{Q},\frac{2m+1}{2}\right)
+\sum_{k=-n}^n\operatorname{Res}\left(\frac{h}{Q},k\right) \right]
$$
implies that
$$
\pi^2 = \sum_{k=-\infty}^{+\infty} 
\frac{1}{\big(k-\frac{2m+1}{2}\big)^2} > \sum_{k=0}^{+\infty} 
\frac{1}{\big(k-\frac{2m+1}{2}\big)^2},
$$
as claimed. This concludes the proof.
\qed

\ms


Finally, we show that the complete system of Theorem
\ref{complete-set} is not a Schauder basis for $A^2(\cW'_\mu)$, for 
all $\mu\ge\pi/2$.    For the definition of Schauder basis, see 
Section~\ref{sec:orthogonal}.  

\begin{thm}\label{red-thm}  \sl
 Let $\mu\ge\pi/2$,  and let   
$H_{\ell,j}(z_1,z_2) =
E_{c_0+\nu\ell+i(j+1)/2}(z_1,z_2)z_2^j$,
$\ell, j\in\bbZ,\,\ell\ge0$.  
For each $j\in\bbZ$ fixed, the function $H_{0,j}$ is in the
$A^2(\cW'_\mu)$-closure of $\operatorname{span} \{ H_{\ell,j},\, \ell=1,2,\dots\}$.   
In particular, this violates the uniqueness requirement in the
definition of Schauder basis.  
\end{thm}

\proof
We first assume that $\mu>\pi/2$.

Let $Q=Q_n$ be a polynomial of degree $n$ of  one complex variable,
without constant term, $Q(w)=\sum_{\ell=1}^n c_\ell w^\ell$.
Then, arguing as in \eqref{F-G} we have
\begin{align} \notag
& \| H_{0,j} - \sum_{\ell=1}^n c_\ell H_{\ell,j} \|_{A^2(\cW'_\mu)}^2
\notag \\
& = \pi \int_{|s|<\mu} \int_\Delta
\big|
\zeta^{c_0+i(j+1)/2} e^{is(c_0+i(j+1)/2)} 
- \sum_{\ell=1}^n c_\ell \zeta^{c_0+\nu\ell+i(j+1)/2}
e^{is(c_0+i(j+1)/2+\nu\ell)}
\big|^2 
\,dA(\z)\,   e^{s(j+1)}  ds \notag \\ 
& = \pi \int_{|s|<\mu} \int_\Delta 
\big| \zeta^{c_0+i(j+1)/2}\big|^2 
\Big|1-\sum_{\ell=1}^n c_\ell   \zeta^{\nu\ell}e^{is\nu\ell}\Big|^2
\,dA(\z)\,    ds \notag \\ 
&\le C  \int_{|s|<\mu} \int_\Delta
\Big| 1 -  \sum_{\ell=1}^n c_\ell   (\z^\nu e^{is\nu})^\ell \Big|^2 
\,dA(\z)\,    ds \notag \\
& = C 
\frac{1}{\nu^2}  \int_{|s|<\mu} \int_{\Omega_s} 
\big| 1 -  \sum_{\ell=1}^n c_\ell  w^\ell \big|^2 
\big| w^{\frac1\nu -1}\big|^2 
\,dA(w)\,     ds \, , \label{RHS-above}
\end{align}
where we have set $w = \z^\nu e^{is\nu}$.
Since  $\mu>\pi/2$, $\nu<1-2\delta$ for some $\delta >0$, 
so that  
$-\frac\pi2< \nu s<\frac\pi2 $ and $0< \nu\frac\pi2 <  \frac\pi2
(1-2\delta)$.  Hence, 
\begin{align*}
\Omega_s 
& \subseteq
\big\{ w = \rho e^{it}:\, 0<\rho<2^\nu,\, 
\nu \big( \textstyle{s-\frac\pi2}\big) <t< 
\nu \big( \textstyle{s+\frac\pi2}\big) 
\big\} \\
& \subseteq \big\{ w = \rho e^{it}:\, 0<\rho<2^\nu,\, 
   |t|< \pi(1-\delta)  
\big\} \\
& =: S\, .
\end{align*}

Plugging this into \eqref{RHS-above} we obtain that
\begin{align*} 
\| H_{0,j} - \sum_{\ell=1}^n c_\ell H_{\ell,j} \|_{A^2(\cW'_\mu)}^2 
& \le C \int_{|s|<\mu} \int_S
\big| 1 -  Q_n(w) \big|^2 
| w|^{2(\frac1\nu -1)} \,dA(w)\,   ds  \\
& =  C \int_S
\big| 1 -  Q_n(w) \big|^2 
| w|^{2(\frac1\nu -1)} \,dA(w)\, .
\end{align*}

Setting $d\omega(w)=| w|^{2(\frac1\nu -1)} dA(w)$,
the conclusion will follow if we show that there exist polynomials
$P_n=1-Q_n$ such that $P_n(0)=1$ and $\|P_n\|_{A^2(S,d\omega)}\to 0$
as $n\to+\infty$.  \ms

In order to prove that such polynomials exist, let $\Delta_+$ be the
half disk $\{z\in\bbC:\, |z|<1,\, \Re z>0\}$ and 
$p(z) = (z-\frac12)^2+\frac34$.   Then $p(0)=1$, and $|p(z)|\le 1$ for
$z\in \overline{\Delta_+}$, as it is elementary to check.  Therefore,
$F(w)=p((2^{-\nu}w)^{ 1/[2(1-\delta)] })$ is a function holomorphic on $S$ such that
\begin{itemize}
\item[{\tiny $\bullet$}] $F$ is continuous on $\overline{S}$;
\item[{\tiny $\bullet$}] $F(0)=1$;
\item[{\tiny $\bullet$}] $|F(w)|\le 1$ on  $\overline{S}$.
\end{itemize}

Observe that $\omega(S)<+\infty$.
Given $\eps>0$, let $K$ be a compact subset of $S$ such that 
$\omega(S\setminus K)<\eps$ and let $n$ be a positive integer such
that $|F^n(w)|\le\eps$ for $w\in K$.  Then
$$
\int_S |F^n(w)|^2 \, d\omega(w) \le C\eps \, .
$$
By Mergelyan's approximation theorem (see \cite{RUD} e.g.), 
we can find polynomials $p_n$
such that $|F^n(w)-p_n(w)|\le \eps$ for $w\in \overline{S}$. Finally,
we set $P_n= \frac{1}{p_n(0)}p_n$ and the conclusion follows easily.
\ms

Finally, let $\mu=\pi/2$, so that $\nu=1$.   Set $\cD=\cup_{|s|<\pi/2}\{z_1:\,
|z_1-e^{is}|<1\}$. 
We have,
\begin{align}
 \| H_{0,j} - \sum_{\ell=1}^n c_\ell H_{\ell,j} \|_{A^2(\cW'_\mu)}^2 
		&= 
\pi \int_{|s|<\mu} \int_{|z_1-e^{is}|<1}
		|z_1^{2c_0+i(j+1)}|\,\big| 1-\sum_{\ell=1}^n c_\ell z_1^{\ell}
                \big|^2\,dA(z_1)\,e^{s(j+1)}\, ds  \notag \\
& 		\le
C \int_\cD
		\big| 1-\sum_{\ell=1}^n c_\ell z_1^{\ell}
                \big|^2\,dA(z_1)\, . \label{Brennan-est}
\end{align}

We observe that $\cD$ is a Jordan domain, 
having  the origin as a boundary point.  
By \cite{Far} we know that the
polynomials are dense in $A^2(\cD)$ and by \cite{Bre85} it follows
that there exists no bounded boundary evaluation point on the space of
polynomials.  Hence, the right hand side of \eqref{Brennan-est} can by
made arbitrarily small and the conclusion now follows.  
 We leave the simple
details to the reader.
\qed

\begin{remark}{\rm
It follows from Theorems \ref{complete-set} 
and \ref{red-thm} 
that the set $\{ H
_{\ell,j}\}_{\ell, j\in\bbZ,\,\ell\ge0}$
is complete, but not a Schauder basis.   It would be of interest to
show that the set $\{ H
_{\ell,j}\}_{\ell, j\in\bbZ,\,\ell\ge0}$ is however a {\em frame} for
$A^2(\cW')$, that is, there exist constants $c_1,c_2>0$ such that 
$$
c_1 \| f\|_{A^2(\cW')}^2 \le 
\sum_{\ell, j\in\bbZ,\,\ell\ge0} \big| \la f, H_{\ell,j}
\ra_{A^2(\cW')} \big|^2 \le c_2 \| f\|_{A^2(\cW')}^2 \,.  \ms
$$
Indeed, the theory of frames in Hilbert function spaces constitute a
fundamental tool, especially in sampling and reconstruction of
functions -- see \cite{DS} where frames were introduced 
in the context of nonharmonic Fourier series,   and  \cite{Young} for
applications 
of the theory of frames to the present setting.   We also recall that
a in a separable Hilbert space, a frame that is also a basis is called
a {\em Riesz basis}.  Hence, in particular, the complete set $\{ H
_{\ell,j}\}_{\ell, j\in\bbZ,\,\ell\ge0}$ of Theorem \ref{complete-set}
is not a Riesz basis either.  
\ms
}
\end{remark}

\section*{Concluding Remarks}

\ms

Thanks to work of several authors, the worm domain has become an important
object of study.  In particular, we are beginning to understand the
Bergman kernel and projection on some versions of the worm.  But the
original smooth worm $\cW_\mu$ is particularly resistive to analysis.
It does not have the built-in symmetries of some of the non-smooth
worms.  In particular, we do not have a useful complete orthogonal
basis for the Bergman space on $\cW_\mu$.  In addition to the
alternative approach mentioned in \cite[\S 5]{KPS1}, this
paper has offered some first steps towards addressing that problem.

As an additional remark, we point out that the results obtained in this work 
can be generalized to the case of worm domains in $\bbC^n$ defined and 
studied in \cite{Barrett-Sahutoglu}.
\ms

\bibliography{worm6-mathscinet}
\bibliographystyle{alpha}

\end{document}